\documentclass[12pt]{article}
\usepackage[english]{babel}
\usepackage{amsmath,amsthm}
\usepackage{amsfonts}
\theoremstyle{definition}

\theoremstyle{remark}

\numberwithin{equation}{section}

\begin{document}

\begin{center}
\bfseries  A COROLLARY  OF  RIEMANN  HYPOTHESIS  \\
\end{center}

\begin{center}

\vspace{5mm}
JINHUA  FEI\\
\vspace{5mm}

ChangLing  Company of Electronic Technology   \, Baoji \,  Shannxi  \,  P.R.China \\
\vspace{5mm}

E-mail: feijinhuayoujian@msn.com \\

\end{center}

\vspace{3mm}

{\bfseries Abstract.}\,  This paper use the  results  of  the  value distribution theory , got  a significant conclusion  by  Riemann hypothesis \\


{\bfseries MR(2000)  Subject  Classification \quad 30D35 , 11M06 } \\

\vspace{8mm}

  First, we give some signs , definition and theorem  in  the value  distribution theory , its contents see the references [1] and [2] . \\

  Definition .

\[  \log^+ x  =   \left\{
     \begin{array}{ccc}
     \log x  \qquad  \quad  1 \leq  x \\
     \quad 0    \qquad  \quad    0 \leq x < 1
     \end{array}
     \right.
\]

It is easy to see that$  \,\,  \log x  \, \leq  \, \log^+ x   $ .\\

  Set $ f(z)  $ is a  meromorphic function  in the region $ \,\, |z|  \, <  \, R \, ,  \,\,  0 \,  < R  \,  \leq  \, \infty  \,\, $  , and not identical to zero .

  $ n(r,f) $  represents  the  poles number  of  $ f(z) $ on  the  circle  $ |z| \, \leq  \, r \,( \, 0 <  r < R \, ) \,  $  ,  multiple  poles  being  repeated . $ n(0,f) $  represents the order of pole of $ f(z)  $  in the  origin .  For arbitrary complex number  $ a  \neq  \infty  \, , \,\,\,\,  n( r , \frac{1}{ f-a  }) \,\, $ represents  the  zeros number  of  $ f(z) - a $ in the circle $ \, |z| \leq  r \,\, (  \, 0 <   r  <  R \, ) \,\,  $ , multiple  zeros  being  repeated. $ \,\,  n( 0 , \frac{1}{f-a} ) \,\, $   represents the order of zero of $ \,\, f(z) - a \,\, $ in the  origin .\\

 Definition .
$$  m( r, f)  =  \frac{1}{ 2 \pi}  \, \int_0^{2 \pi} \, \log^+  \left| f( r e^{i \varphi} )  \, \right|  d \varphi  $$

$$ N ( r , f )\, = \, \int_0^r \, \frac{ n(t,f) - n(0,f)}{ t } \, d t   \,\, + \,\, n(0,f) \, \log \, r $$ \\

 Definition . $ \,\, T( r, f  ) \, = \, m( r , f) \, + \, N ( r , f ) $ .

   $\,\,  T ( r, f ) $  is called the characteristic function  of  $ f ( z )  $ . \\

{\bfseries\footnotesize LEMMA 1.}  If  $ f(z)  $  is  an analytical function in the region  $ \, |z| < R \,\, (\, 0 < R \leq  \infty \, ) $ , then

$$ T ( r , f ) \, \leq \,  \log^+ \, M ( r , f ) \, \leq \,  \frac{ \rho + r }{ \rho - r } \, \, T ( \rho \, , f)    (\, 0 < r < \rho < R  )   $$

where  $ \,\,  M( r, f ) \, = \, \max_{ |z| = r   }  \, | f(z)  | $

The proof of the lemma  see  the  page 57  of the  references [1] . \\

{\bfseries\footnotesize LEMMA 2.}  Set $ f(z) $  is a meromorphic function  in the region $ |z| < R \,\, ( \, 0 < R \leq \infty  )\, $  , not identical to zero .
Set $ |z| < \rho  \,\, ( \, 0 < \rho < R \,) \,\,  $ is a  circle ,  $ \,\, a_\lambda \,\, ( \, \lambda = 1,2, ..., h \, )   $ and  $\,\, b_\mu \,\,  (\, \mu = 1,2,..., k\, ) \,  $ respectively  is the zeros  and  the poles  of    $ f(z)  $  in
   the  circle ,  appeared number of every zero or every pole and its order the same , and that $ z = 0  $ is not the zero  or the pole  of  function   $ f(z)  $  , then in  the   circle  $ |z| < \rho  $ , We have the following formula

$$  \log \, \left| \,  f(0) \,  \right| \, = \,  \frac{1}{ 2 \pi}  \,\int_0^{ 2 \pi } \,\, \log\, \left| f( \rho e^{ i \varphi} )
   \right| d \varphi  \,\, - \,\, \sum_{ \lambda =1}^h \, \log\, \frac{\rho}{|a_\lambda|} \,\, + \,\,
    \sum_{ \mu =1}^k \, \log\, \frac{\rho}{|b_\mu|} $$
this formula is called Jensen formula .

The proof of the lemma  see  the  page 48  of the  references [1] . \\

{\bfseries\footnotesize LEMMA 3.} Set function $ f(z) $ is the meromorphic function  in  $ \, |z| \leq  R \, $  , and

 $$ f(0) \,  \neq  \,\,  0 ,\, \, \infty ,\, \, 1, \,\,\,\, \,\, f'(0) \, \neq \, 0 $$

then when $ \, 0 < r < R \, $ , have

$$  T ( r , f ) \, <  \, 2 \left\{  N ( R, \frac{1}{f}) \, + \, N ( R , f ) \, + \, N ( R , \frac{1}{ f - 1} ) \right\} $$

$$  + \,\, 4\, \, \log^+ |f(0)| \,\, + \,\, 2 \, \, \log^+ \, \frac{1}{ R |f'(0)| } \,\, + \,\, 24\,\, \log \, \frac{R}{R-r}
\,\, + \,\, 2328 $$

This is  a form  of  Nevanlinna second basic theorems .

The proof of the lemma  see the theorem 3.1 of the page 75 of  the  references [1] .  \\

The need for behind, We will make some preparations.\\

{\bfseries\footnotesize LEMMA 4.}  If when $ \, x \geq  a \,  $  ,  $  f(x)$  is  a  nonnegative  degressive  function , then below limits exist

  $$  \lim_{ N \, \rightarrow  \, \infty} \,\,  \left( \, \, \sum_{n=a}^N \, f(n) \,\, - \,\, \int_a^N  f(x) \, d x  \right) \,\,
  = \,\, \alpha $$

where $  \, 0 \leq  \alpha  \leq  f(a)  \, $ . in addition , if when  $ \, x \rightarrow  \infty \,  $ , have $ \, f(x) \rightarrow  0  \, $ , then

$$  \left| \,\,  \sum_{a \leq n \leq \xi }  \, f(n) \,\, - \,\,  \int_a^\xi  f(\nu) \, d \nu   \,\, - \,\, \alpha
 \,\,\right| \,\, \leq \, f( \xi - 1 ) \, , \qquad  (  \, \xi \geq a+1  \,  )  $$

The proof of the lemma  see the theorem 2 of  page 91 of  the  references [3] . \\

    Set $ \, s =  \sigma  +  it \,  $ is the  complex number ,  when $ \, \sigma  > 1  \,  $ , the  definition  of   Riemann  Zeta  function  is

$$  \zeta (s) \, = \, \sum_{n=1}^\infty  \,\, \frac{1}{n^s}   $$

When $ \,  \sigma  >  1 \,  $ , from  the page 90 of   the  references [4], have

$$  \log \zeta (s) \, =  \, \sum_{n=2}^\infty \, \frac{\Lambda(n)}{n^s  \log n}   $$

where $ \, \Lambda (n) \,  $  is  Mangoldt  function . \\

{\bfseries\footnotesize LEMMA 5.}  For any real number  $ t $ , have \\

(1) $$ 0.0426  \, \leq \, \left| \, \, \log \zeta ( 4+ it ) \, \, \right|  \, \leq \, 0.0824   $$

(2)$$  \left| \,\,  \zeta ( 4 + it  )   \, - \, 1  \,\, \right|  \, \geq \, 0.0426   $$

(3)$$  0.917 \, \leq \, \left| \,\,  \zeta  ( 4 + it ) \,\, \right|  \, \leq  \,  1.0824 $$

(4)$$  \left|  \,\,  \zeta ' ( 4 + it ) \,\,   \right|   \, \geq \, 0.012  $$ \\

{\bfseries\footnotesize  PROOF.}

(1) $$    \left| \, \, \log \zeta ( 4+ it ) \,\, \right|  \, \leq \, \sum_{n=2}^\infty \,  \frac{\Lambda(n)}{n^4  \log n} \, \leq  \,
 \sum_{n=2}^\infty \, \frac{1}{n^4} \, = \, \frac{\pi^4}{90} \, - \, 1\, \leq \, 0.0824  $$

$$   \left| \, \, \log \zeta ( 4+ it ) \,\, \right| \, \geq \, \frac{1}{2^4} \, - \, \sum_{n=3}^\infty \, \frac{1}{n^4} \, =
 \, 1 \,+\, \frac{2}{2^4} \, - \,  \sum_{n=1}^\infty \, \frac{1}{n^4} \, = \, \frac{9}{8} \,- \, \frac{\pi^4}{90} \, \geq \, 0.0426  $$

(2)

$$    \left| \, \,  \zeta ( 4+ it )  - 1 \,\, \right| \, = \, \left| \,\, \sum_{n=2}^\infty \, \frac{1}{n^{4 + it}} \,\, \right|
 \, \geq \, \frac{1}{2^4} \, - \, \sum_{n=3}^\infty  \, \frac{1}{n^4}  \,$$

$$= \, 1+ \frac{2}{2^4} \, - \,  \sum_{n=1}^\infty  \, \frac{1}{n^4} \, = \, \frac{9}{8} \, - \, \frac{\pi^4}{90} \, \geq \, 0.0426   $$

(3)

$$   \left| \, \, \zeta ( 4+ it ) \,\, \right| \, = \, \left| \,\sum_{n=1}^\infty  \, \frac{1}{n^{4+it}} \, \right| \, \leq \,
 \sum_{n=1}^\infty  \, \frac{1}{n^4} \, = \, \frac{\pi^4}{90} \, \leq \, 1.0824  $$

$$  \left| \, \, \zeta ( 4+ it ) \,\, \right| \, = \, \left| \,\sum_{n=1}^\infty  \, \frac{1}{n^{4+it}} \, \right| \, \geq  \, 1 \,
 - \,  \sum_{n=2}^\infty  \, \frac{1}{n^4} \, = \, 2 \, - \, \sum_{n=1}^\infty  \, \frac{1}{n^4} \, = \, 2\, - \, \frac{\pi^4}{90} \,
  \geq \, 0.917  $$

(4)

$$    \left| \, \, \zeta '( 4+ it ) \,\, \right| \, =  \, \left| \,\, \sum_{n=2}^\infty \, \frac{\log n }{n^{4+it}}
\,\,\right| \, \geq \, \frac{\log2}{2^4} \, - \,  \sum_{n=3}^\infty \, \frac{\log n }{n^{4}}$$

from lemma 4 , have

$$  \sum_{n=3}^\infty \, \frac{\log n }{n^{4}} \, = \, \int_3^\infty \, \frac{\log x}{x^4} \, d x \, + \, \alpha  $$

where $ \,  0 \, \leq \,  \alpha  \, \leq \, \frac{\log3}{3^4} $

 $$   \int_3^\infty \, \frac{\log x}{x^4} \, d\, x \, = \, - \, \frac{1}{3} \, \int_3^\infty \, \log x \,\, d \, x^{-3} \,
 = \, \frac{\log3}{3^4} \, + \, \frac{1}{3}\, \int_3^\infty \, x^{-4} \, d \, x \, $$

 $$ = \, \frac{\log3}{3^4} \, - \, \frac{1}{3^2} \, \int_3^\infty \, d  \, x^{-3} \, = \, \frac{\log3}{3^4} \, + \, \frac{1}{3^5} $$

therefore

 $$     \sum_{n=3}^\infty \, \frac{\log n }{n^{4}} \, \leq \, \frac{\log3}{3^4} \, + \, \frac{1}{3^5} \, + \, \frac{\log3}{3^4}  $$

therefore

 $$       \left| \, \, \zeta '( 4+ it ) \,\, \right| \, \geq \, \frac{\log2}{2^4} \, - \, \frac{2\log3}{3^4} \, -
 \, \frac{1}{3^5} \, \geq \, 0.012    $$

The proof is complete .\\

   Set $ 0 <  \delta  \leq  \, \frac{1}{100} \, $ , $ \,\,  c_1 ,\, c_2 , \, ... \, , \,  $ represents positive  constant with only $ \delta $  relevant in the article below . \\

{\bfseries\footnotesize LEMMA 6.}  When  $ \,\,  \sigma  \, \geq  \,  \frac{1}{2} , \,\,  |t| \geq  2  \,\, $  , have

$$  \left| \, \zeta ( \sigma  + it  )   \right|  \, \leq  \, c_1  \, |t|^\frac{1}{2}  $$

The proof of the lemma  see the theorem 2 of  page 140 and   the theorem 4  of  page 142 ,  of  the  references [4] . \\

{\bfseries\footnotesize LEMMA 7.} Set $ \, f(z) \, $ is  the  analytic function  in the circle $ \,  | z - z_0  |\, \leq \, R  \,\,$ , then  for any   $ \, 0 < r < R \,   $ , in the circle  $ \,\,   | z - z_0  |\, \leq \, r  \,\, $ , have

$$  \left| \,  f(z) -  f(z_0)  \,  \right|  \, \leq  \,  \frac{2r}{R-r} \, \left ( \, A(R) - Re f(z_0) \,   \right)  $$

where  $ \,\,  A(R)\, = \, \max _{ | z - z_0  | \leq  R }   \,\,  Re f(z)  $

The proof of the lemma  see the theorem 2 of  page 61  of  the  references [4] . \\

   Now assume Riemann hypothesis is correct, abbreviation for RH . In other words , when $ \sigma >  \frac{1}{2}  $ , the  function  $ \zeta ( \sigma  +  i t )  $  has no zeros .  Set  the  union set  of  the region  $ \,  \sigma  >  \frac{1}{2} \, ,  \, | t | > 1  \, $  and  the  region  $ \, \sigma  >  2 \, , \, | t | \leq 1 \, $ is the  region D .

  Therefore , the function $  \zeta ( \sigma + it )  $  have  neither zero nor poles  in the  region D , so ,  function  $  \log  \zeta  ( \sigma + it )  $ is a defined multi-valued analytic function in  the  region D . Every single value analytic branch differ $ 2 \pi i $ integer times .

  Assuming there are the points  $\, s_0 \,$  in the  region D ,  satisfy $ \, \zeta ( s_0 ) = 1 \, $ ( If there is not such point  $\, s_0 \,$ , then the result of   lemma 9 turns into $\,\, N ( \rho , \frac{1}{\zeta-1})\, = 0\,\,$ ,  the  results of the theorem of this article can be obtained directly ). For different single value analytic branch , the  value  of $ \, \log \zeta (s_0) = \log 1 \,  $  are different , it can value  $ \, 0 , \, 2 \pi k i ,( k = \pm 1 , \pm 2, ...... ) \, $ . We select the single valued analytic  branch   of    $  \, \log \zeta (s_0) = \log 1 \, = \, 0   $ .

  Because the region D is simple connected region , so the according to the single value  theorem  of analytic continuation ( the  theorem  see the theorem 2 of  page 276  of  the  references [5] and  theorem 1 of  page 155  of  the  references [6] ),   $  \log  \zeta  ( \sigma + it ) $ is the single valued analytic function  in the region D . in addition ,   when and only  when  $  \zeta  ( \sigma + it )  = 1 $  ,  have $  \log  \zeta  ( \sigma + it ) = 0  $ . In other words ,  1  value   point of  $   \zeta  ( \sigma + it) $  is  the  zero  of  $  \log  \zeta  ( \sigma + it ) $ ,  the  opposite  is true .

  Below, $ \, \log  \zeta  ( \sigma + it ) \,  $   always express a  single valued analytic  branch for we selected .\\

{\bfseries\footnotesize LEMMA 8.} If RH is  correct  ,  then when $ 0 < \delta \leq \frac{1}{100} $  ,  $ \, \sigma\, \geq \, \frac{1}{2} \, + 2 \delta   \, ,\, \, | t |  \geq  16 $ , we have

   $$ \left | \,\, \log  \zeta ( \sigma  + it  ) \,\,  \right | \,\, \leq   \,\,  c_2    \log  | t |  +  c_3   $$

{\bfseries\footnotesize proof.} In the lemma 7  , we choose  $ z_0  = 0 , \,\, f(z) = \log \zeta ( z + 4 + it  )  , \,\, |t| \geq  16  ,  \,\, R = \frac{7}{2}
 - \delta  ,\,\, r = \frac{7}{2}  - 2 \delta  \,\, $ . Because  $ \log \zeta ( z + 4 + it )  $  is the  analytic function  in the circle  $ | z - z_0 | \leq  R  $ , so , from  the lemma 7 ,   in the circle   $ | z - z_0 | \leq  r   $  , we have

$$  \left| \,\,  \log \zeta ( z + 4 + it  ) \,\, - \log \zeta (  4 + it  ) \,\,   \right| \,\, \leq \,\, \frac{7}{ \delta }\,\,
\left (  \,\, A(R)  - Re \log \zeta ( 4 + it ) \,\, \right)  $$

hence

$$   \left| \,\,  \log \zeta ( z + 4 + it  ) \,\, \right|  \leq  \,\, \frac{7}{ \delta }  \,\, \left( \,\,  A(R)  +
 \left| \,\,  \log \zeta ( 4 + it )  \,\,  \right| \,\, \right)  +  \left| \,\, \log \zeta ( 4 + it ) \,\,  \right|  $$

from the lemma 6 , have

$$ A(R) = \max_{  | z - z_0 | \leq  R }   \log \left| \,\, \zeta ( z + 4 + it  ) \,\,  \right|  \leq \,\, \frac{1}{2} \, \log |t|
 \, + \, \log  c_1 $$

from  the lemma 5 , have

$$  \left| \,\,  \log \zeta ( z + 4 + it  ) \,\, \right| \,\, \leq  \,\, c_2  \, \log|t| \, +  \,  c_3  $$

because  $ |t| \geq  16  $  is  real number  arbitrarily , so when  $ \sigma  \geq  \frac{1}{2} + 2 \delta   $ , we have

$$   \left| \,\,  \log \zeta ( \sigma + it  ) \,\, \right| \,\, \leq  \,\, c_2  \, \log|t| \, +  \,  c_3   $$

The proof is complete . \\

{\bfseries\footnotesize LEMMA 9.} If  RH  is  correct , then when $ 0 < \delta \leq \frac{1}{100}  $ ,  $ \, |t| \geq  16 \, ,\, \rho = \frac{7}{2}  -  2 \delta $ ,  in the circle  $ |z| \leq  \rho   $ , we have

$$  N \left( \,\,\rho \, ,\,\, \frac{1}{ \zeta ( z+ 4 + it  )\,\, - 1 } \,\, \right ) \,\, \leq \,\, \log\log |t| \, + \,  c_4  $$ \\

{\bfseries\footnotesize proof.} In the lemma 2 ,  we choose $ f(z) =  \,  \log \zeta ( z+ 4 + it  ) , \,\, R = \frac{7}{2} - \delta ,\,\, \rho = \frac{7}{2} - 2 \delta , \,\,a_\lambda \,\,\, (  \lambda = 1,2, ..., h ) $  is the zeros  of  function  $ \log \zeta ( z+ 4 + it  )   $  in the circle $  | z | < \rho $ , multiple  zeros  being  repeated. The function $ \log \zeta ( z+ 4 + it  ) $   has no poles in the the circle  $  | z | < \rho $ ,  and  $ \log \zeta (  4 + it  )   $  not equal to zero , therefore we have

$$ \log  \left| \,\log  \,  \zeta (  4 + it  ) \, \right|  = \, \frac{1}{ 2 \pi} \, \int_0^{ 2 \pi} \,  \log  \,
\left| \, \log\, \zeta (  4 + it + \rho e^{i \varphi} ) \, \right| \, d\varphi \, - \,  \sum_{\lambda = 1}^{h} \, \log
\, \frac{ \rho }{ |a_\lambda|  } $$

from the lemma 5 and the lemma 8 , have

$$  \sum_{\lambda = 1}^{h} \, \log  \, \frac{ \rho }{ |a_\lambda|  } \,\, \leq \,\,  \log\log |t |  \, + \, c_4  $$

because  $ z = 0  $   is neither the zero ,  nor pole  of  the function $   \log \zeta ( z+ 4 + it  )   $ , so if  $ r_0 $ is a sufficiently small   positive number , then

$$    \sum_{\lambda = 1}^{h} \, \log  \, \frac{ \rho }{ |a_\lambda| } \, = \, \int_{r_0}^\rho  \,
\left( \log \frac{\rho}{t} \right) \,\, d\, n(t, \frac{1}{f}  )  \,\, = \,\, \left[ \left( \log \frac{\rho}{t} \right) \,
   n(t, \frac{1}{f} ) \right] \bigg |_{r_0}^\rho $$

$$ + \,\, \int_{r_0}^\rho  \frac{n(t, \frac{1}{f}  )}{t } \,\, d \, t   \,\, = \,\,\int_{0}^\rho
 \frac{n(t, \frac{1}{f}  )}{t } \,\, d \, t  \,\, = \,\, N \left( \, \rho \, , \frac{1}{f} \,  \right )  $$

$$= \,\, N \left( \, \rho \, , \,\, \frac{1}{\log \zeta ( z+ 4 + it ) }  \, \right )  \,\, = \,\,  N \left( \, \rho \, , \,\,
\frac{1}{ \zeta ( z+ 4 + it ) - 1  }  \, \right )  $$

The proof is complete . \\

{\bfseries\footnotesize THEOREM .} If  RH  is  correct , then  when $  \sigma  \, \geq \, \frac{1}{2}  \,  +  \, 4 \delta \, ,\,\, 0 < \delta \leq
\frac{1}{100}\,, \,\,|t| \geq  16    \,\,  $ , we have

$$  \left| \,  \zeta ( \sigma + it  ) \, \right| \,  \leq    c_{8}  \,\left( \, \log | t | \, \right)^{c_{6}}  $$ \\

{\bfseries\footnotesize proof.}  In the lemma 3 , we choose $ f(z) \,= \, \zeta ( z + 4 + it )  , \, |t| \geq 16  $ , from  the  lemma 5 , have  $\, f(0) \, = \, \zeta ( 4 + it ) \, \neq \,  0, \, \infty , \,  1, \,\,\,\, \,\,  f'(0) \, = \,  \zeta ' (  4 + it ) \, \neq  0   $ , and  $   f'(0) \, = \,  \zeta ' (  4 + it
) \, \geq  0.012  \, , \,\,\,\,  | \, f(0) \, |  \, = \, \left| \, \zeta ( 4 + it ) \, \right|  \, \leq  \, 1.0824 $ . We choose  $ R \, = \,  \frac{7}{2} \, - \, 2 \delta , \,\, r \, = \, \frac{7}{2} \,  - \, 3 \delta  $ .  because $ \zeta ( z + 4 +it )  $  is  the analytic function , and have neither  zero nor the poles in the circle $ |z| \leq  R $ , therefore

$$   N \left( \, R \, , \, \frac{1}{f}  \right)  \,  =\, 0 \,\, , \,\,\,\,\,\,\,  \,\,\,\,  N  \left( \, R \, , \, f  \right )  =  0    $$

from the lemma 9 , have

$$  T \, \left( \, r \, , \zeta ( z + 4 + it  )  \, \right ) \,\,  \leq  \,\,  2 \log\log |t| \, + \, c_5    $$

In the lemma 1 ,  we choose $ \, R = \frac{7}{2} - 2\delta \,, \,  \rho \, = \, \frac{7}{2}  \, - \,  3 \delta  , \,\,  r \, = \,  \frac{7}{2} \, -
\, 4 \delta   $ , from the maximal principle , in the the circle  $ |z | \leq  r   $ , we have

$$  \log^+  \left| \,\zeta ( z + 4 + it  ) \,  \right|   \, \leq  \, c_6 \, \log\log |t| \, + \, c_7   $$

Since $  |t|\geq 16  $ is arbitrary real  number, so when  $ \sigma  \geq \frac{1}{2} \, + \, 4 \delta   $ , have

$$   \log^+  \left| \,\zeta ( \sigma + it  ) \,  \right|   \, \leq  \, c_6 \, \log\log |t| \, + \, c_7   $$

therefore

$$  \log \,  \left| \,\zeta ( \sigma + it  ) \,  \right|   \, \leq  \, c_6 \, \log\log |t| \, + \, c_7     $$

therefore

$$  \left| \,  \zeta ( \sigma + it  ) \, \right| \,  \leq    c_{8}  \,\left( \, \log | t | \, \right)^{c_{6}}  $$

The proof is complete . \\

The result of this theorem  is  better than   known   results .

\vspace{20mm} \centerline{ REFERENCES } \vspace{5mm}

[1]  Zhuang Q.T,  Singular   direction   of   meromorphic   function ,  BeiJing: Science Press,1982 . ( in Chinese )   .\\

[2]  Yang  L  ,  Value   distribution   theory    and   new   research ,    BeiJing: Science Press,1982. ( in Chinese ) \\

[3]  Hua  L.G  ,   Introduction  of  number  theory  ,   BeiJing: Science Press,1979. ( in Chinese ) \\

[4]  Pan C.D, Pan C.B, Fundamentals of analytic number theory,    BeiJing:  Science Press, 1999 . ( in Chinese )  \\

[5]  Zhuang  Q.T , Zhang  N.Y ,   Complex variables functions , BeiJing:  Peking University press ,  1984 . ( in Chinese )  \\

[6]  Hua  L.G ,    Introduction  of   advanced   mathematics  ( Book  One  of second volume  ) ,  BeiJing:  Science Press, 1981 . ( in Chinese )

\end{document}